\documentclass[12pt]{article}
\usepackage{amsmath,amssymb,amsthm}
\usepackage[utf8]{inputenc}
\usepackage{fontenc}
\RequirePackage[top=1in,bottom=1in,left=1in,right=1in,includehead]{geometry}
\usepackage{aliascnt}
\usepackage{hyperref,cleveref}
\usepackage{tikz}
\usepackage{genyoungtabtikz}

\newcommand{\CSn}{\mathbb{C}S_n}
\DeclareMathOperator{\eig}{eig}
\DeclareMathOperator{\diag}{diag}

\DeclareMathOperator{\type}{type}
\DeclareMathOperator{\shape}{shape}
\DeclareMathOperator{\noninv}{noninv}

\newcommand{\tab}{\Yboxdim{12pt}\young}

\newcommand{\Nu}[1]{\nu_{#1}}
\newcommand{\allnuk}{\{\Nu{k}\}_{k \in \mathbb{N}}}

\newtheorem{thm}{Theorem}
\newaliascnt{lem}{thm}
\newaliascnt{ex}{thm}
\newaliascnt{prop}{thm}
\newaliascnt{corl}{thm}
\newaliascnt{rmq}{thm}

\newtheorem{ex}[ex]{Example}
\newtheorem{prop}[prop]{Proposition}
\newtheorem{corl}[corl]{Corollary}
\newtheorem{rmq}[rmq]{Remark}
\aliascntresetthe{lem}
\aliascntresetthe{ex}
\aliascntresetthe{prop}
\aliascntresetthe{corl}
\aliascntresetthe{rmq}

\numberwithin{equation}{section}

\title{Eigenvalues of symmetrized shuffling operators}

\author{Nadia Lafrenière 
\thanks{\href{mailto:lafreniere.nadia.2@courrier.uqam.ca}{lafreniere.nadia.2@courrier.uqam.ca}.
The author was supported by a graduate scholarship from Fonds de recherche Québec and Institut des Sciences mathématiques.}\addressmark{1}}

\def\@address{\relax}
\def \address{\@getaddress}
\def \@getaddress#1{{
        \gdef \@address{#1}
    }}
    
\newcommand \addressmark[1]{%
    $^{\ #1}$%
}

\address{\addressmark{1}Laboratoire de combinatoire et d'informatique mathématique, Université du Québec à Montréal, Québec, Canada}

\date{\today}

\begin{document}

\maketitle
{\small \itshape \@address \par}

\textbf{Abstract:} This paper describes a combinatorial way of obtaining all the eigenvalues of the symmetrized shuffling operators introduced by Victor Reiner, Franco Saliola and Volkmar Welker. It allows us to prove their conjecture that these eigenvalues are integers. This work generalizes the case of the random-to-random Markov chain.

\textbf{Résumé :} Cet article décrit une stratégie combinatoire pour calculer toutes les valeurs propres des opérateurs de mélange symétrisés introduits par Victor Reiner, Franco Saliola et Volkmar Welker. Ceci nous permet de prouver leur conjecture affirmant que ces valeurs propres sont entières. Les travaux présentés ici généralisent le cas du mélange doublement aléatoire.


\section{Context}

The \textit{random-to-random shuffle} on a sequence of objects $(w_1, w_2, \ldots, w_n)$ can be described as follows: pick any object randomly, with uniform probability, remove it, and place it back anywhere in the sequence once again randomly. Introduced by Persi Diaconis and Laurent Saloff-Coste \cite{DSC}, it has drawn attention in the literature on shufflings ever since \cite{DS,QM,SCZ,Subag,UR}, in part because the computation of its mixing time proved elusive. The mixing time is the number of times one must execute a shuffle to get a sequence that is well shuffled. It was finally computed by Megan Bernstein and Evita Nestoridi \cite{BN}. The key ingredient in their proof was a combinatorial description of the eigenvalues of the random-to-random Markov chain found by Anton Dieker and Franco Saliola \cite{DS} using representation theory and the following interpretation of the Markov chain. The \textit{random-to-top shuffle} consists of taking any object in our sequence (randomly) and reinserting it at the beginning of the sequence. The \textit{top-to-random shuffle} is the dual version of random-to-top and is done by taking the first element in the sequence and inserting it in the sequence at a random position. The random-to-random shuffle is simply the composition of random-to-top and then top-to-random. It is said to be the \textit{symmetrized version} of random-to-top.

One can generalize the random-to-random shuffle by replacing the random-to-top operator by another one moving more than one object at a time. Victor Reiner, Franco Saliola and Volkmar Welker observed such a family of operators $\allnuk$ that have special properties \cite{RSW}. These operators $\nu_k = \pi_k^\top\pi_k^{}$ are such that $\pi_k$ follows the move-and-randomize scheme, which means that one needs to remove $k$ objects randomly, and put them on top in a random order. Although the eigenvalues of the $\{\pi_k\}_{k\in\mathbb{N}}$ have been known for a long time \cite{BHR, DFP}, this does not give the eigenvalues of the $\allnuk$. The challenge lies in the fact that very few properties of a Markov chain are shared with its symmetrized version.

As said earlier, the operators in $\allnuk$ are of a certain interest not only because they generalize the random-to-random shuffle, but also for what Reiner, Saliola and Welker found about them: they commute, and they only have real non-negative eigenvalues. They also conjectured that all their eigenvalues are integers (up to a rescaling), a statement we confirm in \autoref{corl:evals_intergers}. Furthermore, we show how to explicitly compute all the eigenvalues of those operators in \autoref{thm:main}.

This result has the potential to give good bounds on the mixing time of the Markov chains associated to $\allnuk$. Some techniques to bound the mixing time use the eigenvalues \cite{DiSh, LPW}. 

As well, the results presented in this paper give a general framework to study the eigenvalues of all the symmetrized shuffling operators described by Reiner, Saliola and Welker \cite{RSW}. In their paper, they described another infinite family of operators $\{\gamma_k\}_{k\in \mathbb{N}}$ that pairwise commute, and they conjectured that no other pair of symmetrized shuffling operators commute. The techniques we show here could be used to give a combinatorial description for the eigenvalues of  $\{\gamma_k\}_{k\in \mathbb{N}}$. 

The paper is organized as follows. \Cref{sec:desc_op} describes the operators and \autoref{sec:eval} presents the necessary background to compute the eigenvalues. The latter section also gives the formula for computing the eigenvalues. In \autoref{sec:content}, we explore what happens if multiple occurrences of the same object appear in the sequence. Finally, \autoref{sec:theory} describes the connections with representation theory. The reader who is curious about the connections between tableaux and the eigenvalues of the shuffling operators is invited to jump ahead to \autoref{sec:theory} before reading the paper.

\section{Description of the operators}\label{sec:desc_op}
Let $w = (w_1,\ldots, w_n)$ be a sequence of objects. The operator $\Nu{k}$ removes $k$ objects from $w$ and reinserts them one after the other, so that it does not necessarily preserve the order of those $k$ objects. We encode the result of that operation as the formal linear combination of all possible permutations of $w_1, \ldots, w_n$ one can obtain doing this shuffle. In that linear combination, the coefficient of a shuffled sequence is the number of ways of obtaining it after one shuffle. 
\begin{ex}
	Applying $\Nu{1}$ to the sequence \textsf{abc} means removing one of the three letters and reinserting it at one of the three remaining spaces. For example, if one removes the letter \textsf{a}, then \textsf{abc}, 
	\textsf{bac} and \textsf{bca} can be obtained in one manner each. Doing the same with \textsf{b} and \textsf{c}, one obtains
	\[\Nu{1}(\textsf{abc}) = 3 \cdot \textsf{abc} + 2 \cdot \textsf{bac} + 2 
	\cdot \textsf{acb} + \textsf{cab} + \textsf{bca}.\]
\end{ex}

\begin{ex}
	The operator $\Nu{2}$ acts on the sequence \textsf{abcd} by moving two letters, with letters that can be moved to the same place. The coefficient of \textsf{cabd} in $\Nu{2}(\textsf{abcd})$ is $4$. The pairs of letters to be moved to obtain \textsf{cabd} are either $\{\textsf{a,b}\}$, $\{\textsf{a,c}\}$, $\{\textsf{b,c}\}$ or $\{\textsf{c,d}\}$.
\end{ex}
The operator $\Nu{0}$ is the identity, since there is only one way to remove no object in the sequence.

\begin{rmq}
	In \cite{RSW}, Reiner, Saliola and Welker denote these operators $\nu_{(n-k,1^k)}$. In their context, it is natural to use this notation since they look at symmetrized shuffling operators $\{\nu_{\lambda}\}_{\lambda \vdash n}$ indexed by partitions. Here, we only use partitions of the form $(n-k, 1^k)$, and we always specify the value of $n$, if needed.
\end{rmq}

The definition of those operators in terms of shuffling is handy for developing some intuition, but it is hard to write explicitly the matrices of the linear operators from this description. To solve that problem, one can use the following definition, taken from \cite{RSW}. Let $\sigma \in S_n$ be a permutation and let $\noninv_k(\sigma)$ be the number of increasing sequences of length $k$ in $\sigma$. Then, $\Nu{k}$ acts on $\sigma$ by 
\[ \Nu{k}(\sigma) = \sum_{\tau \in S_n} \noninv_{n-k}(\sigma^{-1}\tau) \tau. \]

\begin{ex}\label{ex:mat_Nu1}
	The matrix for the linear operator $\Nu{1}$ on sequences of length 3 is
	\[\bordermatrix{
		&{\scriptstyle \textsf{abc}} & {\scriptstyle \textsf{acb}} & 
		{\scriptstyle \textsf{bac}} & {\scriptstyle \textsf{bca}} & 
		{\scriptstyle \textsf{cab}} & {\scriptstyle \textsf{cba}} \cr
		{\scriptstyle \textsf{abc}} & 3 & 2 & 2 & 1 & 1 & 0 \cr
		{\scriptstyle \textsf{acb}} & 2 & 3 & 1 & 2 & 0 & 1 \cr
		{\scriptstyle \textsf{bac}} & 2 & 1 & 3 & 0 & 2 & 1 \cr
		{\scriptstyle \textsf{bca}} & 1 & 2 & 0 & 3 & 1 & 2 \cr
		{\scriptstyle \textsf{cab}} & 1 & 0 & 2 & 1 & 3 & 2 \cr
		{\scriptstyle \textsf{cba}} & 0 & 1 & 1 & 2 & 2 & 3 \cr
	}.\]
\end{ex}

Note that the matrix in \autoref{ex:mat_Nu1} is not a transition matrix. Since the entries in a transition matrix are transition probabilities, they always take values between 0 and 1, and all rows sum to 1. 
From the operators $\allnuk$, one can recover a transition matrix by dividing  the entries of the matrix of $\nu_k$ on sequences of length $n$ by  $\binom{n}{k}^2 k!$. One can verify that the matrix in \autoref{ex:mat_Nu1} has all rows summing to $9$. This scaling of the matrix is what allows us to have integer eigenvalues. Recall that the complex eigenvalues of a transition matrix of a Markov chain all have modulus smaller than or equal to 1.

\section{Eigenvalues}\label{sec:eval}
The eigenvalues of the symmetrized shuffling operator $\Nu{k}$ are obtained recursively, using the eigenvalues of $\Nu{k-1}$. Fortunately, the eigenvalues for $\Nu{1}$, the random-to-random shuffle, are given by a formula presented in \autoref{sec:eval_R2R}, found by Dieker and Saliola \cite{DS}. This formula is given in terms of statistics on tableaux. Hence, the first part of this section is dedicated to presenting the combinatorial notions necessary to understand the computations.

\subsection{Partitions, diagrams and tableaux}
Given a positive integer $n$, a \textit{partition $\lambda$ of $n$} is a decreasing list of positive integers $(\lambda_1, \ldots, \lambda_r)$ whose sum is $n$. To each partition $\lambda = (\lambda_1, \ldots, \lambda_r) \vdash n$, one can associate a \textit{diagram} formed by $\lambda_1$ left-aligned boxes in the first row, $\lambda_2$ left-aligned boxes in the second row, and so on, the rows being placed from top to bottom, like in matrices. Examples of partitions and diagrams are found in \autoref{ex:diagram}. 

Given two diagrams $\lambda$ and $\mu$ such that the boxes of $\mu$ are contained in $\lambda$, the \textit{skew diagram} $\lambda/\mu$ is the set of boxes of $\lambda$ that do not belong to $\mu$. As pictured in \autoref{ex:diagram}, it is not necessarily left-justified. 

Given a diagram, one can write positive integers in the boxes and thus obtain a \textit{tableau}. A tableau is said to be \textit{standard} if the numbers in the $n$ boxes are the integers from 1 up to $n$ and if they are placed in the following way: from left to right and from top to bottom, the numbers in each row and column are strictly increasing. All standard tableaux of size $2$ and $3$ are listed in \autoref{ex:SYT}. 

\begin{ex}\label{ex:diagram}
	The tuple $\lambda = (4,3,2,1)$ is a partition of $10$ and can be drawn as in the first picture. The second picture shows a tableau of shape $(3,1,1)$, while the third one is the skew diagram $(4,3,2,1)/(3,1,1)$.
	\[ \Yboxdim{16pt} \yng(4,3,2,1)\ , \qquad \young(933,1,4)\ ,\qquad 
	\gyoung(:::;,:;;,:;,;) \ .\]
\end{ex}

\begin{ex}\label{ex:SYT}
	The standard tableaux of size 2 and 3 are
	\[\Yboxdim{16pt}\young(12)\ ,\quad \young(1,2)\ , \quad \young(123)\ ,
	\quad \young(12,3)\ , \quad \young(13,2)\ , \quad  \young(1,2,3)\ .\]
\end{ex}

A statistic on a (skew) diagram that is useful to our computation is its \textit{diagonal index}. To each cell in a diagram, one can associate a number by taking its coordinates $(i,j)$, like in a matrix, and computing the difference $j-i$ (see \autoref{ex:diagonal_index}). For a (skew) diagram $\lambda$, its diagonal index is the sum of the diagonal indices of its cells and is denoted $\diag(\lambda)$. 
\begin{ex} \label{ex:diagonal_index}
	The diagonal indices of the cells in the following diagrams are written in each cell. The left diagram has diagonal index $7$, while the right one has diagonal index $17$.
	\[\Yboxdim{16pt}
	\yngres(0,5,4,4,1)\ , \qquad  \yngres(0,6,4)\ .\]
\end{ex}
\begin{rmq}
	Some texts refer to the diagonal index as the \textit{content} of a cell. However, the content of a tableau can also qualify the numbers that fill the cells. To remove ambiguity, we only use the word content with the latter meaning (defined more explicitly in \autoref{sec:content}).
\end{rmq}

\subsection{Schützenberger's $\Delta$ operator}
The strategy we use to compute the eigenvalues is an induction both on the operator's index, $k$, and on the number of objects in our sequences, $n$. The idea is that it is not hard to compute the eigenvalues of the matrices for very small values of $n$. But since the dimensions of the matrices grow factorially, it is too much to expect that we can compute the eigenvalues naively. To do such an induction on $n$, we need a way to find, from a tableau of size $n$, the right tableau of size $n-1$ to proceed with the induction. 
This is what the \textit{Schützenberger $\Delta$ operator} does 
\cite{schutzenberger}. To execute it on a tableau $t$ of size $n$:
\begin{enumerate}
	\item We remove the element 1 and replace it with an empty box.
	\item While the empty box is not in an outer corner of the tableau, we move it toward the outside using jeu-de-taquin slides: this means that if $(i,j)$ is the empty box, we switch it with the box containing the minimum of the boxes at position $(i+1,j)$ and $(i,j+1)$.
	\item Once the empty box is in an outer corner, we get a new tableau without the hole that contains the numbers $2, \ldots, n$. In order to recover a standard tableau, we replace those numbers with $1, \ldots, n-1$, preserving the order, and remove the empty box.
\end{enumerate}
An example of this process appears in \autoref{ex:jdt}.

\begin{ex}\label{ex:jdt}
	Applying the Schützenberger $\Delta$ operator to the standard tableau $t$ of size $7$ returns a standard tableau of size $6$.
	\[t = \Yboxdim{16pt}\young(125,347,6) \to 
	\young(\bullet25,347,6) \to \young(2\bullet5,347,6) \to 
	\young(245,3\bullet7,6) \to \young(245,37\bullet,6) \to 
	\young(134,26,5) = \Delta(t). \]
\end{ex}

This process allows one to deduce a unique standard tableau of size $n-1$ from a standard tableau of size $n$.

\subsection{Horizontal strips and desarrangement tableaux}
A (skew) diagram is a \textit{horizontal strip} if no two cells lie in the same column.

In a standard tableau of size $n$, an entry $i$ is an \textit{ascent} if either $i = n$ or if the box containing $i+1$ is weakly to the north-east of the one containing $i$.
A standard tableau is a \textit{desarrangement tableau} if it is empty or if 
its smallest ascent is even. Many equivalent definitions exist as well, see for instance \cite{DW, RSW}.

\begin{ex}\label{ex:horiz_strip}
	The skew diagram $(4,2,2,1)/(3,2,1)$ is a horizontal strip, as shown in this diagram on the left, where the cells from the strip are identified by a cross. 
	The ascents in the right tableau are identified by a shaded background.
	\newcommand{\ylw}{\Yfillcolour{yellow}}
	\newcommand{\wh}{\Yfillcolour{white}}
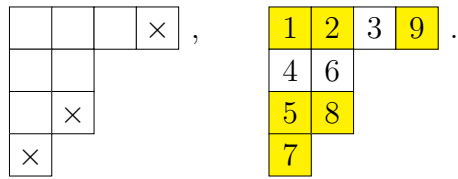
	\[\Yboxdim{16pt}\gyoung(;;;;\times,;;,;;\times,\times)\ , \qquad
	\young(!\ylw12!\wh3!\ylw9,!\wh46,!\ylw58,7)\ .\]
\end{ex}
\begin{ex}\label{ex:desar_tableaux}
	Among the tableaux from \autoref{ex:SYT}, only the second and the fifth are desarrangement tableaux.
\end{ex}

\begin{prop}[Proposition 6.25 in \cite{RSW}]\label{prop:desa_vs_bandes}
	To each standard tableau $t$, one can associate a unique desarrangement tableau by applying the Schützenberger $\Delta$ operator until the result is a desarrangement tableau. If one needs to do it $j$ times, then $t/\Delta^j(t)$ is a horizontal strip of size $j$.
\end{prop}

Using \autoref{prop:desa_vs_bandes}, one can define the \textit{type} of a standard tableau $t$ as the minimum number $j$ such that $\Delta^j(t)$ is a desarrangement tableau.

\begin{ex}
	The tableau of \autoref{ex:horiz_strip} has type $3$ and 
	\[\Delta^3(t) =  \Yboxdim{16pt}\young(136,25,4)\ .\]
	The diagram of the shape of $t/\Delta^3(t)$ is the horizontal strip from 
	\autoref{ex:horiz_strip}.
\end{ex}

\subsection{Eigenvalues of the random-to-random operator} \label{sec:eval_R2R}
Anton Dieker and Franco Saliola recently described the eigenvalues of the  random-to-random shuffle \cite{DS}. This shuffling operator is not only a very natural shuffle and a special case of the operators we are describing (it is a rescaling of $\Nu{1}$), it also serves as the base case of our induction for the eigenvalues of $\allnuk$. It is useful to recall their description of the eigenvalues:
\begin{thm}[Theorem 5 in \cite{DS}]
	Every eigenvalue of the random-to-random shuffle is of the form $\frac{1}{n^2}\eig(\lambda/\mu)$, where $\lambda$ is a partition of $n$ and $\lambda/\mu$ is a horizontal strip. Moreover, $\eig$ is the following combinatorial statistic defined on skew partitions:
	\[\eig(\lambda/\mu) = \binom{|\lambda| + 1}{2} - \binom{|\mu| + 1}{2} 
	+ \diag(\lambda/\mu). \]
	
	The multiplicity of an eigenvalue $\frac{1}{n^2}\varepsilon$ is
	\[ \sum_{\substack{ \lambda/\mu\text{ is a horizontal strip,}\\
			\eig(\lambda/\mu)=\varepsilon}} f^\lambda d^\mu,\]
	where $f^\lambda$ is the number of standard tableaux of shape $\lambda$, and $d^\mu$ is the number of desarrangement tableaux of shape $\mu$.
\end{thm}
\begin{rmq}
	Using this theorem, the reader might notice that all the eigenvalues are real and between $0$ and $1$. In that context, $\frac{1}{n^2}\eig(\lambda/\mu)$ is not an integer, but rather the eigenvalue of the random-to-random Markov chain. The integer eigenvalue of $\Nu{1}$ is given by $\eig(\lambda/\mu)$.
\end{rmq}
As well, one can associate a standard tableau with an eigenvalue of the random-to-random shuffle \cite{Franco}. This reformulation will be used in \autoref{sec:main_result} to compute the eigenvalues for other operators. For a standard tableau $t$:
\begin{itemize}
	\item If $t$ is a desarrangement tableau, the eigenvalue associated to $t$ is 0.
	\item Otherwise, the eigenvalue is $\varepsilon + n + \diag(t/\Delta(t))$, where $\varepsilon$ is the eigenvalue associated to $\Delta(t)$.
\end{itemize}
Since any standard tableau leads to a desarrangement tableau through multiple iterations of $\Delta$ (by \autoref{prop:desa_vs_bandes}), this process always terminates.

\subsection{Main result: Eigenvalues of $\Nu{k}$}\label{sec:main_result}

We are now ready to present our main result. It describes all the eigenvalues of the operators $\allnuk$. Using a result from Reiner, Saliola and Welker, the kernel of $\nu_k$ is associated with the tableaux of type smaller than $k$ (see \autoref{eq:kernel_nu}, or Theorem VI.9.5 in \cite{RSW}). Thus, to all tableaux of type at least $k$, we must associate a non-zero eigenvalue. The way we compute this value allows us to show it is a positive integer (\autoref{corl:evals_intergers}).

\begin{thm}\label{thm:main}
	The eigenvalues of $\Nu{k}$ on $\CSn$ are indexed by the standard tableaux of size $n$.
	For such a tableau $t$, the eigenvalue $v_k(t)$ is given by
	\[v_k(t) = \begin{cases}
	v_k(\Delta(t)) + \left(n+1-k+\diag(t/\Delta(t))\right)\ 
	v_{k-1}(\Delta(t)) & \text{ if } \type(t) \geq k,\\
	0 & otherwise.
	\end{cases}\]
	The tableau $t$ contributes $f^{\shape(t)}$ to the multiplicity of the eigenvalue $v_k(t)$, where $f^{\shape(t)}$ is the number of standard tableaux of the shape of $t$.
\end{thm}

Note that, replacing $k$ by $1$, one can recover the eigenvalues for the  random-to-random operator, presented in \autoref{sec:eval_R2R}.

\begin{ex}\label{ex:main}
	In this example, we compute all the eigenvalues of $\Nu{2}$ as a random walk on $S_4$. Since the Markov chain has a state space of size $4! = 24$, one can still compute its eigenvalues with the usual algorithms. We did so using Sagemath \cite{sagemath}, and we therefore know that the eigenvalues are $0$, $4$, $20$ and $72$, with multiplicities $17$, $3$, $3$ and $1$, respectively. These values are confirmed by \autoref{thm:main}.\\
	There are ten standard tableaux of size $4$. 
	\[\Yboxdim{14pt}
	{\setlength\arraycolsep{3.5pt}
		\renewcommand{\arraystretch}{1.3}
		\begin{array}{|p{16mm}|c|c|c|c|c|c|c|c|c|c|}
		\hline
		\text{Standard} \newline \text{tableau} & \young(1234) 
		& \young(123,4) & 
		\young(124,3) & 
		\young(134,2) & 
		\young(12,34) & \young(13,24) & \young(12,3,4) & \young(13,2,4) & 
		\young(14,2,3) 
		& \young(1,2,3,4)\\
		\hline
		\textrm{Type} & 4 & 2 & 1 & 0 & 1 & 0 & 2 & 0 & 1 & 0\\
		\hline
		\end{array}}\]
	\Yvcentermath1  
	To find the non-zero eigenvalues of $\Nu{2}$, one needs only to compute $v_k(t)$ for the standard tableaux of type at least 2. There are 3 such tableaux:
	\allowdisplaybreaks
	\begin{align*}
	v_2\left(\tab(1234)\right)
	&= v_2\left(\tab(123)\right) + (4+1-2+3)\cdot 
	v_1\left(\tab(123)\right)\\
	&= \left(  v_2\left(\tab(12)\right) + (3+1-2+2)\cdot 
	v_1\left(\tab(12)\right) \right) + 6 \cdot 
	v_1\left(\tab(123)\right)\\
	&= \left(\left(v_2\left(\tab(1)\right) + (2+1-2+1) 
	\cdot v_1\left(\tab(1)\right) \right) + 4\cdot 
	v_1\left(\tab(12)\right) 
	\right) + 6 \cdot 9\\
	& = \left(\left( 0 + 2 \cdot v_1\left(\tab(1)\right)\right)+ 4 
	\cdot 4 
	\right) + 54\\
	& = 2 \cdot 1 + 16+ 54 = 72.\\
	v_2\left(\tab(123,4)\right)
	&= v_2\left(\tab(12,3)\right) + (4+1-2+2)\cdot 
	v_1\left(\tab(12,3)\right)\\
	&= 0 + 5 \cdot v_1\left(\tab(12,3)\right)\\
	&= 0 + 5 \cdot 4 = 20.\\
	v_2\left(\tab(12,3,4)\right) &= v_2\left(\tab(1,2,3)\right) 
	+(4+1-2+1) \cdot v_1\left(\tab(1,2,3)\right)\\
	& = 0 + 4 \cdot v_1\left(\tab(1,2,3)\right)\\
	& = 4 \cdot 1 = 4.
	\end{align*}

	In the last computation, the second equality comes from the fact that $\tab(1,2,3)$ has type smaller than 2, while the third one has been computed using the result from \autoref{sec:eval_R2R}.
\end{ex}

The work presented here allows us to confirm Conjecture I.1.2, found in \cite{RSW}.
\begin{corl}\label{corl:evals_intergers}
	All eigenvalues of the operators $\allnuk$ are (positive) integers.
\end{corl}

\subsection{Sequences with multiple copies of the same object}\label{sec:content}  
Up to now, we viewed our sequences as lists containing objects that are all distinct. However, those operators also work on sequences with repetitions. For example, it is common to apply it to a word, and words often have repeated letters. The Markov chains are defined in the same way, but they have fewer states. For example, the matrix for $\Nu{2}$ on a sequence of four 
elements, two of each of the two kinds, is
\[ \Nu{2} = \bordermatrix{
	&{\scriptstyle \textsf{aabb}} & {\scriptstyle \textsf{abab}} & {\scriptstyle \textsf{baab}} & {\scriptstyle \textsf{abba}} &
	{\scriptstyle \textsf{baba}} & {\scriptstyle \textsf{bbaa}}\cr
	{\scriptstyle \textsf{aabb}} & 20 & 16 & 12 & 12 & 8 & 4 \cr
	{\scriptstyle \textsf{abab}} &16 & 14 & 12 & 12 & 10 & 8 \cr
	{\scriptstyle \textsf{baab}} & 12 & 12 & 12 & 12 & 12 & 12 \cr
	{\scriptstyle \textsf{abba}} & 12 & 12 & 12 & 12 & 12 & 12 \cr
	{\scriptstyle \textsf{baba}} & 8 & 10 & 12 & 12 & 14 & 16 \cr
	{\scriptstyle \textsf{bbaa}} & 4 & 8 & 12 & 12 & 16 & 20\\
}.\]

What we call the \textit{content} of a sequence is the quantity of each of the objects in the sequence, placed in decreasing order. It is a partition, and this allows us to compare it to other partitions.

If $\lambda = (\lambda_1, \ldots, \lambda_l)$ and $\mu = (\mu_1, \ldots,  \mu_m)$ are two partitions of $n$, we say $\lambda$ \textit{dominates} $\mu$, written $\lambda \unrhd \mu$, if, for all $i \in \{1, 2, \ldots, \min(l,m)\}$,
\[ \lambda_1 + \ldots + \lambda_i \geq \mu_1 + \ldots + \mu_i. \]

\begin{thm}\label{thm:eval_content}
	The non-zero eigenvalues of $\Nu{k}$ on words of content $\mu \vdash n$ are indexed by the standard tableaux of size $n$, type at least $k$ and shape $\lambda \unrhd \mu$. For such a tableau $t$,
	\[v_k(t) = v_k(\Delta(t)) + (n+1-k+\diag(t/\Delta(t)))\ 
	v_{k-1}(\Delta(t)).\]
\end{thm}
This agrees with \autoref{thm:main}, since permutations have content $(1,1,\ldots,1)$, which is minimal for the dominance order.

\begin{ex}
	The matrix for $\Nu{2}$ on words of content $(2,2)$, like $\textsf{aabb}$, has eigenvalues $72$, $20$ and $0$ with multiplicities $1$, $1$ and $4$, respectively. The content of the words being $(2,2)$, we need to find the partitions that dominate it. Those are $(4)$, $(3,1)$ and $(2,2)$. To get the eigenvalues, we use \autoref{thm:eval_content} and what we computed in \autoref{ex:main}.
\end{ex}

\section{How it works: The connection between tableaux and eigenvalues of shuffling operators}\label{sec:theory}
This section is dedicated to explaining briefly the reason we can use tableaux to compute the eigenvalues. Tableaux are frequently used in the representation theory of the symmetric group.

If we consider the case where there is no repetition in the sequence, then the operators act on the symmetric group algebra. It is well-known that the simple modules of the symmetric group are the Specht modules $S^\lambda$, indexed by partitions. This way, one can decompose $\CSn$:
\[ \CSn \cong \bigoplus_{\lambda\vdash n} f^\lambda S^\lambda \cong \bigoplus_{t \text{ a standard tableau}} S^{\shape(t)},\]
where $f^\lambda$ is the number of standard tableaux of shape $\lambda$. Hence, there is a bijection between copies of simple modules and standard tableaux, and one can associate each standard tableau with a copy of a Specht module.  For more explanations, see for example \cite{sagan}.

Moreover, since Specht modules are simple, Schur's lemma shows that their morphisms (shufflings, for example) have a single eigenvalue when restricted to a Specht module. Combined with the remark in the last paragraph, we associate each tableau with an eigenvalue.

Another useful decomposition of $\CSn$ is the decomposition into the \textit{permutation modules} $\{M^\lambda\}_{\lambda \vdash n}$, indexed by partitions:
\[\CSn \cong \bigoplus_{\lambda \vdash n} M^\lambda.\]
Those modules are not simple, and thus may carry more than one eigenvalue for a given shuffling operator, yet they have some utility. The basis of $M^\lambda$ is the set of words of content $\lambda$, i.e.\ the words containing $\lambda_1$ occurrences of the first letter, $\lambda_2$ occurrences of the second, and so on. This means we can directly work on words to understand the action of the operators $\nu_k$. To understand the impact of the Schützenberger $\Delta$ operator, we compare the effect of $\Nu{k}$ on $M^\lambda$ and $M^{\lambda-e_i}$ (when the only cell in $t/\Delta(t)$ is in the $i$-th row). Then, taking the restriction to the Specht module $S^{\lambda-e_i}$, we get an exact computation of how the eigenvalues evolve in that process. Finally, this led to \autoref{thm:main} for the non-zero eigenvalues.

To distinguish when a simple module belongs to the kernel or to another eigenspace for a given operator $\Nu{k}$, we use two observations. The first one concerns specifically the sequences with repeated elements, for example those with content $\mu$. The reason we restrict our study to tableaux whose shape dominate $\mu$ is given by Young's rule: this gives the decomposition of the permutation modules into Specht modules. The multiplicity of $S^\lambda$ in $M^\mu$ is $m_{\lambda, \mu}$, the number of semistandard tableaux of shape $\lambda$ and content $\mu$:
\[ M^\mu \cong \bigoplus_{\lambda \vdash n} m_{\lambda, \mu} S^\lambda. \]
Recall that a semistandard tableau is a tableau whose entries are strictly increasing on the columns and weakly increasing on the rows. It is not hard to see that $m_{\lambda,\mu} = 0$ when $\lambda \ntrianglerighteq \mu$. For more details, see for example Chapter 2 in \cite{sagan}.

The second observation we use on the kernel comes from Reiner, Saliola and Welker and their study of $\allnuk$ and their kernels. They showed that the Schützenberger $\Delta$ operator connected tableaux of size $n$ and type $j$ and  tableaux of size $n-1$ and type $j-1$, but also eigenspaces of the symmetrized shuffling operators for sequences of $n$ and $n-1$ elements. The trio of authors also explained the construction of the kernel in terms of tableaux of a given type. They found that the kernel of $\Nu{k}$ is 
\begin{equation}
\ker(\Nu{k}) \cong \bigoplus_{\substack{t \text{ a standard tableau,}\\ \type(t) < k}} S^{\shape(t)}. \label{eq:kernel_nu}
\end{equation}
For more details, see sections VI.9 to VI.11 in \cite{RSW}. 

\newcommand \acknowledgements{\section*{Acknowledgements}}
\acknowledgements{I would like to thank my advisor Franco Saliola and my colleague Stéphanie Schanck for the very helpful discussions on that project and for the initial work. Aram Dermenjian and Pauline Hubert read a previous version of this paper and I am grateful for their suggestions. I am appreciative of the comments the referees provided; they improved the article.}

\bibliography{biblio_doc.bib}

\end{document}